\documentclass{article}
\usepackage{amsmath,amsthm,amsfonts,eucal}

\numberwithin{equation}{section}

\def\ca{{\mathcal A}}
\def\cb{{\mathcal B}}

\def\cd{{\mathcal D}}

\def\cf{{\mathcal F}}

\def\cj{{\mathcal J}}
\def\ck{{\mathcal K}}

\def\cn{{\mathcal N}}

\def\cu{{\mathcal U}}

\def\bc{{\mathbb C}}

\def\bn{{\mathbb N}}
\def\br{{\mathbb R}}

\def\a{\alpha}
\def\b{\beta}
        \def\G{\Gamma}
\def\d{\delta}        \def\D{\Delta}
\def\eps{\varepsilon}

\def\th{\vartheta}

\def\l{\lambda}       \def\La{\Lambda}
\def\m{\mu}

\def\t{\tau}

\def\f{\varphi}

\def\om{\omega}        \def\O{\Omega}

\def\imply{\Rightarrow}

\def\e#1{{\rm e}^{#1}}
\def\itm#1{\item[$(#1)$]}

\DeclareMathOperator{\inj}{inj}
\DeclareMathOperator{\supp}{supp}

\def\dec{\searrow}

\DeclareMathOperator{\Lim}{Lim}

\def\lo{\Lim_{\om}}

\newtheorem{Thm}{Theorem}[section]
\newtheorem{Cor}[Thm]{Corollary}
\newtheorem{Prop}[Thm]{Proposition}
\newtheorem{Lemma}[Thm]{Lemma}
\theoremstyle{definition}
\newtheorem{Dfn}[Thm]{Definition}

\theoremstyle{remark}
\newtheorem{rem}[Thm]{Remark} 
 
 \title{\huge A semicontinuous trace for almost local operators on 
 an open manifold}
 \author{Daniele Guido$^{1}$, Tommaso Isola$^{2}$\\
 $(1)$ Dipartimento di Matematica,\\ Universit\`a della Basilicata,\\ 
 I--85100 Potenza, Italy.\\
 $(2)$ Dipartimento di Matematica,\\ Universit\`a di Roma ``Tor Vergata'',\\ 
 I--00133 Roma, Italy.\\
 {\tt guido@unibas.it,\ isola@mat.uniroma2.it}}
\date{April 23, 2001}
\begin{document}
\maketitle
\markboth{A semicontinuous trace package}
{A semicontinuous trace on an open manifold}
\renewcommand{\sectionmark}[1]{}
\bigskip

\begin{abstract} 
	A semicontinuous semifinite trace is constructed on the 
	C$^{*}$-algebra generated by the finite propagation operators 
	acting on the L$^{2}$-sections of a hermitian vector bundle on an 
	amenable open manifold of bounded geometry.  This trace is the 
	semicontinuous regularization of a functional already considered by 
	J.~Roe.  As an application, we show that, by means of this 
	semicontinuous trace, Novikov-Shubin numbers for amenable 
	manifolds can be defined.
\end{abstract}

\newpage
\setcounter{section}{-1}
\section{Introduction.}\label{sec:intro}

 In this paper, we construct a C$^{*}$-algebra of operators on the 
 Hilbert space of $L^{2}$-sections of any vector bundle on an open 
 manifold with bounded geometry, together with a semicontinuous 
 semifinite trace on it which is finite on locally trace-class 
 operators with suitable uniform bounds.  In some sense our 
 C$^{*}$-algebra will be maximal with respect to these properties.

 The finite-dimensionality of the eigenspaces of elliptic differential 
 operators on a compact manifold is one of the classical results 
 allowing us to define of many analytical quantities which then turn 
 out to have a geometric meaning.  The possibility of defining 
 analogous invariants for open manifolds is related to the possibility 
 of renormalizing the dimension of such spaces, in order to give a new 
 sense to finiteness.

 Since the dimension of a space can be defined as the trace of its 
 orthogonal projection, one can try to perform the renormalization by 
 replacing the usual trace with another one.  As eigenprojections of 
 elliptic differential operators on an open manifold are not compact 
 in general, no trace on $B(H)$ can be finite on them, so the required 
 trace can only be defined on a sub-algebra.

 The first step in this direction is due to Atiyah in the seminal 
 paper on the index theorem for covering manifolds \cite{Atiyah}.  He 
 observed that, when the manifold has the structure of an infinite 
 covering of a compact manifold with respect to a group $\G$, 
 $\G$-periodic operators form a type II$_{\infty}$ von Neumann 
 algebra.  If the size of the eigenprojections of $\G$-periodic 
 elliptic differential operators is measured with the normal trace on 
 the von~Neumann algebra of $\G$-periodic operators, the classical 
 finiteness is recovered.  This started the whole theory of 
 L$^{2}$-invariants, producing both new versions of classical 
 invariants, such as L$^{2}$-Betti numbers, and completely new 
 objects, such as Novikov-Shubin numbers.
 
 The idea of extending Atiyah's construction, without any reference to 
 a group, is due to John Roe \cite{Roe} and is based on the following 
 procedure: choose a suitable invasion of the manifold via compact 
 sets (amenable exhaustion), normalize the trace on such sets dividing 
 by the measure of the set, go to infinity by a suitable generalized 
 limit.  Applying this procedure to finite propagation operators 
 gives a trace.

 In order to produce a semicontinuous semifinite trace on a 
 C$^*$-algebra containing finite propagation operators, some technical 
 problems have to be solved, and this is the purpose of this paper.  
 Thus, we choose as a C$^*$-algebra simply the norm closure of the 
 finite propagation operators, in some sense a maximal choice if one 
 requires the trace property.  Further we observe that with the naive 
 definition of the functional is not semicontinuous and has a weak 
 trace property, so we regularize to obtain a semicontinuous 
 semifinite positive trace on the C$^*$-algebra.
 
 
 While the C$^{*}$-algebra only depends on the metric on the 
 manifold $M$, the trace $Tr_{\ck,\om}$ depends upon the large scale 
 geometry of $M$, described by the amenable exhaustion $\ck$, and on a 
 generalized limit procedure $\om$. We expect geometrically 
 interesting quantities derived from  $Tr_{\ck,\om}$ to be independent 
 of $\om$. 
 
 The domain of the trace $\cd(Tr_{\ck,\om})$ contains all compact 
 operators, and indeeed $Tr_{\ck,\om}$ vanishes on them, for any $\om$.  
 Moreover $\cd(Tr_{\ck,\om})$ contains integral operators with 
 off-diagonal exponentially decaying kernels, and, if the kernel is also 
 uniformly continuous on the diagonal, the trace of the integral 
 operator can be computed using the simple procedure devised by 
 J.~Roe.

 As $C_{c}$-functional calculi of $\D_{p}$ belong to 
 $\cd(Tr_{\ck,\om})$, the spectral density function $N_{p}$ of $\D_{p}$ 
 is defined.  Indeed, as $f\in C_{c}[0,\infty) \to 
 Tr_{\ck,\om}(f(\D_{p}))\in\bc$ is a positive linear functional, it 
 defines a Radon measure $\m_{p}$ on 
 $[0,\infty)$, hence we set $N_{p}(t) := \m_{p}((0,t))$. 
 
 As an application, one can define, in analogy with the case of 
 covering manifolds \cite{Atiyah,NS}, both the L$^{2}$-Betti numbers 
 and the Novikov-Shubin numbers of $M$.  Note that L$^{2}$-Betti 
 numbers were also defined (differently but equivalently) by J.~Roe 
 \cite{RoeBetti}.  We prove here that, for open manifolds of bounded 
 geometry, the $0$-th Novikov-Shubin number is always $\geq 1$ thus 
 extending a result of Varopoulos \cite{Varopoulos1}.
 
 As it is known, a general understanding of the geometric meaning of 
 the Novikov-Shubin invariants is still lacking.  We believe that some 
 aspects can be better understood by interpreting them as global 
 invariants of an open manifold, rather than as homotopy invariants of 
 a compact one.
 
 The asymptotic character of these numbers manifests itself at two 
 points in their construction.  On the one hand, the trace used to 
 define these numbers is a large scale trace, since, as observed by 
 Roe \cite{Roe}, it is given by an average over the group, in the case 
 of coverings, and by an average on the exhaustion, in the case of 
 open manifolds.  On the other hand, these numbers are defined in terms 
 of the low frequency behaviour of the $p$-Laplacians, or the large 
 time behaviour of the $p$-heat kernel.  In this respect, they are the 
 large scale counterpart of the spectral dimension, namely of the 
 dimension as defined by the Weyl asymptotics.
 
 On the technical side, the functional $\f$ considered by J.~Roe has 
 two main problems.  The first is that it is not semicontinuous; in 
 fact its kernel is not closed (see Proposition \ref{3.4.1}); the 
 second is that the equality $\f(ab)=\f(ba)$ has only been proved for 
 $a,\ b$ uniformly smoothing operators of order $-\infty$ in the 
 domain of $\f$.  Conversely, the trace property on a C$^{*}$-algebra 
 requires the equality $\f(ab)=\f(ba)$ to hold for $a$ in the domain 
 and $b$ in the C$^{*}$-algebra.  This is proved for the regularized 
 trace in Theorem \ref{def:trace}.
 
 These two properties play a key role in showing that Novikov-Shubin 
 numbers are invariant under quasi-isometries and that they are 
 asymptotic dimensions in the sense of noncommutative geometry, which 
 is done in \cite{GI4}.
 
 Some of the results contained in the present paper have been 
 announced in several international conferences.  In particular we 
 would like to thank the Erwin Schr\" odinger Institute in Vienna, 
 where a preliminary version of this paper was completed, and the 
 organisers of the ``Spectral Geometry Program'' for their kind 
 invitation.

 \section{Open manifolds of bounded geometry}
 \label{subsec:openmanifold}

 \subsection{Preliminaries}
 In this subsection we give some preliminary results on open manifolds
 of bounded geometry that are needed in the sequel.

 Several definitions of bounded geometry for an open manifold (i.e. a
 noncompact complete Riemannian manifold) are usually considered.
 They all require some uniform bound (either from above or from below)
 on some geometric objects, such as: injectivity radius, sectional
 curvature, Ricci curvature, Riemann curvature tensor etc.  (For all
 unexplained notions see e.g. Chavel's book \cite{Chavel}).

 In this paper the following form is used, but see
 \cite{ChavelFeldman} and references therein for a different approach.

 \begin{Dfn}\label{2.1.1}
	 Let $(M,g)$ be a complete Riemannian manifold.
	 We say that $M$ has C$^{\infty}$-bounded geometry if it has
	 positive injectivity radius, and the curvature tensor is bounded,
	 together with all its covariant derivatives.
 \end{Dfn}

 \begin{Lemma}\label{2.1.2}
	 Let $M$ be an $n$-dimensional complete Riemannian manifold with
	 positive injectivity radius, sectional curvature bounded from
	 above, and Ricci curvature bounded from below, in particular M
	 could have C$^\infty$-bounded geometry.  Then there are real
	 functions $\b_1,\ \b_2$ s.t.
	 \itm{i} for all $x\in M$, $r>0$,
	 $$
	 0<\b_1(r)\leq vol(B(x,r)) \leq \b_2(r),
	 $$
	 \itm{ii} $\lim_{r\to0}\frac{\b_2(r)}{\b_1(r)} =1$.
 \end{Lemma}
 \begin{proof}
	 $(i)$ We can assume, without loss of generality, that the
	 sectional curvature is bounded from above by some positive
	 constant $c_1$, and the Ricci curvature is bounded from below by
	 $(n-1)c_2 g$, with $c_{2}<0$.  Then, denoting with $V_\d(r)$ the
	 volume of a ball of radius $r$ in a manifold of constant
	 sectional curvature equal to $\d$, we can set $\b_1(r) :=
	 V_{c_1}(r\wedge r_0)$, and $\b_2:= V_{c_2}(r)$, where $r_0:=
	 \min\{ \inj(M), \frac{\pi}{\sqrt{c_1}} \}$, and $\inj(M)$ is the
	 injectivity radius of $M$.  Then the result follows from
	 (\cite{Chavel}, p.119,123).  \\
	 \noindent $(ii)$
	 \begin{align*}
		 \lim_{r\to0} \frac{\b_2(r)}{\b_1(r)} & = \lim_{r\to0}
		 \frac{V_{c_2}(r)}{V_{c_1}(r)} = \lim_{r\to0} \frac{ \int_0^r
		 S_{c_2}(t)^{n-1}dt }{ \int_0^r S_{c_1}(t)^{n-1}dt } \\
		 & = \left( \lim_{r\to0} \frac{S_{c_2}(r)}{S_{c_1}(r)}
		 \right)^{n-1} = 1
	 \end{align*}
	 where (cfr.  \cite{Chavel}, formulas (2.48), (3.24), (3.25))
	 $V_\d(r) = \frac{n\sqrt{\pi} }{ \G(n/2+1)} \int_0^r
	 S_\d(t)^{n-1}dt$, and
	 $$
	 S_\d(r) :=
	 \begin{cases}
		 \frac1{\sqrt{-\d}}\sinh(r\sqrt{-\d}) & \d<0 \\
		 r& \d=0 \\
		 \frac1{\sqrt{\d}}\sin(r\sqrt{\d}) &  \d>0.
	 \end{cases}
	 $$
 \end{proof}

 \medskip

 Let $M$ be a complete Riemannian manifold, and recall (\cite{Wolf})
 that $\D_{p}:=(d+d^{*})^{2}|_{L^{2}(\La^{p}T^{*}M)}$, the $p$-th
 Laplacian on $M$, is essentially self-adjoint and positive, and the
 semigroup $\e{-t\D_{p}}$ has a C$^\infty$ kernel, $H_{p}(t,x,y)$, on
 $(0,\infty)\times M\times M$, called the $p$-th heat kernel.  Let us
 mention the following result, which will be useful in the sequel.

 \begin{Prop}\label{2.1.5}
	 Let $M$ be an $n$-dimensional complete Riemannian manifold with
	 C$^{\infty}$-bounded geometry, then for all $T>0$, there are $c,
	 c'>0$, s.t., for $0<t\leq T$,
	 \begin{align*}
		 |H_{p}(t,x,y)| & \leq c\ t^{-n/2-1}
		 \exp \left(\frac{-c'\d(x,y)^{2}}{t}\right)  \\
		 |\nabla_{x}H_{p}(t,x,y)| & \leq c\ t^{-n/2-3/2}
		 \exp \left(\frac{-c'\d(x,y)^{2}}{t}\right)
	 \end{align*}
	 where we denoted with $\d$ the metric induced on $M$ by $g$.  As
	 a consequence $H_{p}(t,\cdot,\cdot)$ is uniformly continuous on a
	 neighborhood of the diagonal of $M\times M$.
 \end{Prop}
 \begin{proof}
	 The estimates are proved in \cite{BE}.  For the last statement,
	 for any $\d_{0}<\min\{1,\inj(M),\frac{\pi}{ \sqrt{c_{1}} } \}$,
	 $x\in M$, $y\in B(x,\d_{0})$, we have
	 $|H_{p}(t,x,y)-H_{p}(t,x,x)|\leq \sup|\nabla_{y}H_{p}(t,x,y)|
	 \d(x,y)$, and we get the uniform continuity.
 \end{proof}

 \subsection{The C$^*$-algebra of almost local operators}
 \label{subsec:almostlocal}

 Here we introduce the C$^{*}$-algebra of almost local
 operators, and observe that the $p$-heat semigroup belongs to it.

 Let $F$ be a finite dimensional Hermitian vector bundle over $M$, and
 let $L^{2}(F)$ be the Hilbert space completion of the smooth sections
 with compact support of $F$ w.r.t. the scalar product $\langle
 s_{1},s_{2} \rangle := \int_{M} \langle s_{1x},s_{2x} \rangle
 dvol(x)$. The following lemmawill beused in the sequel.
 
  \begin{Lemma}\label{3.1.4}
	 Let $A$ be a bounded self-adjoint operator on $L^2(F)$, with
	 measurable kernel.  Then
	 $$
	 \|A\| \leq \sup_{x\in M} \int_M |a(x,y)|dy
	 $$
 \end{Lemma}
 \begin{proof}
	 Since $A$ is self-adjoint, $a(x,y)$ is symmetric, hence
	 \begin{align*}
		 \|A\|_{1\to1} & = \sup \{|(f,Ag)| : f\in L^\infty(F),\
		 \|f\|_\infty=1,\ g\in L^1(F),\ \|g\|_1=1 \} \\
		 & \leq \sup_{x\in M} \int_M |a(y,x)|dy =
		 \|A\|_{\infty\to\infty}
	 \end{align*}
	 The thesis easily follows from Riesz-Thorin interpolation
	 theorem.
 \end{proof}

 Recall \cite{RoeCoarse} that an operator $A\in\cb(L^2(F))$ has finite
 propagation if there is a constant $u_A>0$ s.t. for any compact
 subset $K$ of $M$, any $\f\in L^2(F)$, $\supp\f\subset K$, we have
 $\supp A\f \subset Pen^+(K,u_A) := \{ x\in M : \d(x,K)\leq u_A \}$.
 Clearly finite propagation operators for a $^{*}$-subalgebra of 
 $\cb(L^2(F))$, denoted by $\ca_0\equiv \ca_{0}(F)$.  
 
%
%
 \begin{Dfn}
	 We call the norm closure of $\ca_0$ the C$^*$-algebra of almost 
	 local operators on $L^{2}(F)$, and denote it by $\ca\equiv \ca(F)$.
 \end{Dfn}

 \begin{Thm}\label{3.1.6}
	 Let $M$ be a complete Riemannian manifold of C$^{\infty}$-bounded
	 geometry.  Then
	 the C$^{*}$-algebra of almost local operators on
	 $L^{2}(\La^{p}T^{*}M)$ contains all compact operators and the 
	 $C_0([0,\infty))$-functional
	 calculus of the Laplace operator on $p$-forms, i.e.
	 $f(\D_{p})\in\ca(\La^{p}T^{*}M)$, for any $f\in C_0([0,\infty))$
 \end{Thm}
  \begin{proof}
	 Indeed J. Roe shows that compact operators (\cite{RoeCoarse}, 
	 Lemma 4.12) and $C_0([0,\infty))$-functional calculi of $\D_{p}$ 
	 (\cite{RoeCBMS}, Proposition 3.6) belong to the C$^{*}$-algebra 
	 $C^{*}(M)$, which consists of the locally compact operators in 
	 $\ca$.
 \end{proof}

 \section{A functional described by J.~Roe}
 \label{sec:roe}

 This section is devoted to the construction of a trace on the 
 C$^{*}$-algebra $\ca=\ca(F)$.  The basic idea for this construction 
 is due to Roe \cite{Roe}, and is based on a regular exhaustion for 
 the manifold.  We shall regularize this functional, in order to get a 
 semicontinuous semifinite trace on the C$^*$-algebra of almost local 
 operators.  As observed by Roe, this trace is strictly related to the 
 trace constructed by Atiyah \cite{Atiyah} in the case of covering 
 manifolds.  It may therefore be used to define the Novikov-Shubin 
 invariants for open manifolds, as we do in section 
 \ref{sec:regularized}.

 In the rest of this paper $M$ is a complete Riemannian manifold
 of C$^\infty$-bounded geometry as in Definition \ref{2.1.1}, that we
 assume endowed with a regular exhaustion, as in the following 
 Definition.
 
 \begin{Dfn}{\rm \cite{Roe}} 
	 A regular exhaustion $\ck$ of $M$ is an increasing sequence 
	 $\{K_{n}\}$ of compact subsets of $M$, whose union is $M$, and 
	 s.t., for any $r>0$
	 $$									
	 \lim_{n\to\infty} \frac{vol(K_{n}(r)) }{vol(K_{n}(-r))} =1,
	 $$
	 where $K(r)\equiv Pen^{+}(K,r):= \{x\in M: \d(x,K)\leq r\}$, and 
	 $K(-r)\equiv Pen^{-}(K,r):=$ the closure of $M\setminus 
	 Pen^{+}(M\setminus K,r)$.
 \end{Dfn}
 
 Observe that, as $M$ is complete, $Pen^+(K,r)$ coincides with the 
 closure of $\{ x\in M : \d(x,K) < r \}$, which is the original 
 definition of Roe.
  
 \begin{Lemma}\label{l:penombra} 
	 Let $K$ be a compact subset of $M$, then
	 \itm{i} $K(-r_{2})\subset K\subset K(r_{1})$, for any $r_1,r_2>0$
	 \itm{ii} $\{x\in M: \d(x,M\setminus K)<r\} \subset$ Interior of 
 $Pen^+(M\setminus K,r) \equiv M\setminus K(-r)$ 
	 \itm{iii} $Pen^{+}(K(r_{1})\setminus K(-r_{2}),R) \subset 
	 K(r_{1}+R+\eps)\setminus K(-r_{2}-R-\eps)$, for any 
	 $r_1,r_2,R,\eps>0$.
 \end{Lemma}
 \begin{proof}
	 $(ii)$ If $\d(x,M\setminus K)<r$, there is $z\in M\setminus K$ 
	 s.t. $\d(x,z)<r$, so that $x$ belongs to the interior of 
	 $Pen^{+}(M\setminus K,r)$, which is the complement of $K(-r)$.\\
	 $(iii)$ Indeed if $x\in Pen^{+}(K(r_{1})\setminus K(-r_{2}),R)$, 
	 then for any $\eps>0$ there is $x_{\eps}\in K(r_{1})\setminus 
	 K(-r_{2})$ with $\d(x,x_\eps)<R+\eps/2$.  Therefore, on the one 
	 hand, $\d(x,K) \leq R + \eps/2 + r_{1}$, which implies $x\in 
	 K(r_{1}+R+\eps)$.  On the other hand, as $x_{\eps}\not\in 
	 K(-r_{2})$, there is $y_{\eps}\in M\setminus K$ s.t. 
	 $\d(y_{\eps},x_{\eps})<r_{2}+\frac{\eps}2$, hence 
	 $\d(x,y_{\eps})\leq \frac{\eps}2 + R + r_{2}+\frac{\eps}2$ and 
	 $x\not\in K(-r_{2}-R-\eps)$.
 \end{proof}
  
 Following Moore-Schochet \cite{Moore-Schochet}, we recall that an 
 operator $T$ on $L^2(F)$ is called locally trace class if, for any 
 compact set $K\subset M$, $E_KTE_K$ is trace class, where $E_K$ 
 denotes the projection given by the characteristic function of $K$.  
 It is known that the functional $\m_T(K):= Tr(E_KTE_K)$ extends to a 
 Radon measure on $M$.  To state the next definition we need some 
 preliminary notions.

 \begin{Dfn} 
	 Define $\cj_{0+}\equiv \cj_{0+}(F)$ as the set of positive 
	 locally trace class operators $T$, such that 
	 \itm{i} there is $c>0$ s.t. $\m_T(K_{n}) \leq c\ vol(K_{n})$, 
	 for all $n\in\bn$, 
	 \itm{ii} $\lim_{n\to\infty}\frac{\m_{T}(K_{n}(r_{1})\setminus 
	 K_{n}(-r_{2}))}{vol(K_{n})} = 0$, for all $r_{1},r_{2}>0$.
 \end{Dfn}
 
 \begin{Lemma}\label{3.2.3} 
	 $\cj_{0+}$ is a hereditary (positive) cone in $\cb(L^2(F))$.
 \end{Lemma}
 \begin{proof} Linearity follows by $\mu_{A+B}=\mu_A+\mu_B$.
	 If $T\in\cj_{0+}$, and $0\leq A\leq T$, then $Tr(BAB^*)\leq 
	 Tr(BTB^*)$, for any $B\in\cb(L^2(F))$, and the thesis follows. 	 
 \end{proof}
 
 \begin{rem}\label{rem:laplacian} 
	 The hereditary cone $\cj_{0+}$ depends on the exhaustion $\ck$, 
	 however it contains a (hereditary) subcone, given by the 
	 operators $T$ for which there is $c>0$ such that 
	 $\m_T(\Omega)\leq c\ vol(\Omega)$ for any measurable set 
	 $\Omega$.  Proposition~\ref{2.1.5} implies that the operator 
	 $e^{-t\Delta_p}$ belongs to the subcone, hence to 
	 $\cj_{0+}(\La^{p}T^*{M})$.
 \end{rem}

 Recall \cite{Roe} that $\cu_{-\infty}(F)$ is the set of uniform 
 operators of order $-\infty$.
 
 \begin{Prop}\label{prop:cu}
 	$\cu_{-\infty}(F)_{+} \subset \cj_{0+}(F)$.
 \end{Prop}
 \begin{proof}
 	Let $A\in \cu_{-\infty}(F)$, so that $Au(x) = \int_{M} a(x,y)u(y) 
 	dy$, with $a\in C^{\infty}(F\otimes F)$ is a smoothing kernel, and 
 	is uniformly bounded together with all its covariant derivatives 
 	(\cite{Roe}, 2.9). Then for any Borel set $\O\subset M$, 
 	$\m_{A}(\O) = Tr(E_{\O}AE_{\O}) = \int_{\O} tr(a(x,x))dx \leq c\ 
 	vol(\O)$, and the result easily follows.
 \end{proof}
 
 If $\om$ is a state on $\ell^\infty(\bn)$ vanishing on infinitesimal 
 sequences, we use in the following the notation $\lo a_{n}:= 
 \om(\{a_{n}\})$, for any $\{a_{n}\}\in\ell^\infty(\bn)$.  Consider the 
 weight $\f\equiv \f_{\ck,\om}$ on $\cb(L^2(F))_+$ given by
 $$
 \f(A):= 
 \begin{cases}
	 \lo \frac{\m_A(K_{n})}{ vol(K_{n})}  &  A\in\cj_{0+} \\
	 +\infty &  A\in\cb(L^2(F))_+\setminus\cj_{0+}.
 \end{cases}
 $$
 Observe that the functional $\f$ is the functional defined by Roe 
 in \cite{Roe}, but for the domain.
  
 \begin{Prop}\label{p:Roeformula}
	For any $A\in\cu_{-\infty}(F)_+$, $\f(A) = \lo \frac{\int_{K_{n}} 
	tr(a(x,x)) dx}{vol(K_{n})}$, which is Roe's definition in {\rm 
	\cite{Roe}}.
 \end{Prop}
 \begin{proof}
 	Follows easily from the proof of Proposition \ref{prop:cu}.
 \end{proof}
 
 \begin{Lemma}\label{l:phiomega}
	 If $A\in\cj_{0+}$ then 
	 \begin{equation*}
		\f(A) = \lo \frac{\m_A(K_{n}(r_{1}))}{ vol(K_{n}(r_{2}))}
	 \end{equation*}
	 for any $r_1,r_2\in\br$.
 \end{Lemma}
 \begin{proof}
	 Indeed, if $r_{1}\geq0$, we get
	 \begin{equation*}
		 \lo\frac{\m_A(K_{n}(r_{1}))}{ vol(K_{n}(r_{2}))} =\lo\left( 
		 \frac{\m_A(K_{n})}{ vol(K_{n})} + 
		 \frac{\m_A(K_{n}(r_{1})\setminus K_{n})}{ vol(K_{n})}\right) 
		 \frac{vol(K_{n})}{ vol(K_{n}(r_{2}))} = \f(A)
	 \end{equation*}
	 whereas, if $r_{1}<0$, we get
	 \begin{equation*}
		\lo\frac{\m_A(K_{n}(r_{1}))}{ vol(K_{n}(r_{2}))} =\lo\left( 
		\frac{\m_A(K_{n})}{ vol(K_{n})} - \frac{\m_A(K_{n}\setminus 
		K_{n}(r_{1}))}{ vol(K_{n})}\right) \frac{vol(K_{n})}{ 
		vol(K_{n}(r_{2}))} =\f(A).
	 \end{equation*}
 \end{proof}
 
 The algebra $\ca$, being a C$^*$-algebra, contains many unitary 
 operators, and is indeed generated by them.  The algebra $\ca_0$ may 
 not, but all unitaries in $\ca$ may be approximated by elements in 
 $\ca_0$.  Such approximants are $\d$-unitaries, according to the 
 following
 
 \begin{Dfn}\label{3.2.5} 
	 An operator $U\in\cb(L^2(F))$ is called $\d$-unitary, $\d>0$, if 
	 $\|U^*U-1\|<\d$, and $\|UU^*-1\|<\d$.
 \end{Dfn}

 Let us denote with $\cu_\d$ the set of $\d$-unitaries in $\ca_0$ and 
 observe that, if $\d<1$, $\cu_\d$ consists of invertible operators, 
 and $U\in\cu_\d$ implies $U^{-1}\in\cu_{\d/(1-\d)}$.

 \begin{Prop}\label{3.2.6} 
	 The weight $\f$ is $\eps$-invariant for $\d$-unitaries in 
	 $\ca_0$, namely, for any $\eps\in(0,1)$, there is $\d>0$ s.t., 
	 for any $U\in\cu_\d$, and $A\in\ca_+$,
	 $$
	 (1-\eps)\f(A) \leq \f(UAU^*) \leq (1+\eps)\f(A) .
	 $$
 \end{Prop}

 \begin{Lemma}\label{3.2.7} 
	 If $T\in\cj_{0+}$, then $ATA^*\in\cj_{0+}$ for all $A\in\ca_0$.
 \end{Lemma}
 \begin{proof}
	 First observe that for any Borel set $\Omega\subset M$ we have
	 \begin{align*}
		 \m_{ATA^*}(\Omega)
		 & = Tr(E_\Omega ATA^*E_\Omega) \\
		 & = Tr(E_\Omega AE_{\Omega(u_A)}TE_{\Omega(u_A)}A^*E_\Omega) \\  
		 & \leq \|A^*E_\Omega A\| Tr(E_{\Omega(u_A)}TE_{\Omega(u_A)}) \\  
		 & \leq \|A\|^2 \m_{T}(\Omega(u_A))
	 \end{align*}
	 so that 
	 \begin{align*}
		 \frac{\m_{ATA^*}(K_{n})}{vol(K_{n})}
		 & \leq \|A\|^2 \frac{\m_{T}(K_{n}(u_A))}{vol(K_{n})} \\
		 & = \|A\|^2 \frac{\m_{T}(K_{n})}{vol(K_{n})} + 
		 \|A\|^2\frac{\m_{T}(K_{n}(u_A)\setminus K_{n})}{vol(K_{n})}
	 \end{align*}
	 which is bounded.  Now observe that, by Lemma 
	 \ref{l:penombra} $(iii)$, it follows
	 \begin{align*}
		 \frac{\m_{ATA^*}(K_{n}(r_{1})\setminus K_{n}(-r_{2}))}{vol 
		 K_{n}} & \leq \|A\|^{2} 
		 \frac{\m_{T}(Pen^{+}(K_{n}(r_{1})\setminus 
		 K_{n}(-r_{2}),u_A))}{vol K_{n}} \\
		 & \leq \|A\|^{2} \frac{\m_{T}(K_{n}(r_{1}+u_A+\eps)\setminus 
		 K_{n}(-r_{2}-u_A-\eps))}{vol K_{n}} \to 0,
	 \end{align*}
	 the thesis follows.
 \end{proof}
 
 \medskip
  
 \noindent {\it Proof of Proposition \ref{3.2.6}.}
 Assume $A\in \cj_{0+}\cap\ca_{+}$, then $UAU^{*}\in\cj_{0+}$ and, by 
 Lemma \ref{l:phiomega},
 \begin{align*}
	 \f(UAU^*) & = \lo \frac{\m_{UAU^{*}}(K_{n})}{ vol(K_{n})} \\
	 & \leq \|U\|^2 \lo\left(\frac{\m_A(K_{n}(u_{U})) }{ vol(K_{n})} 
	 \right)\\
	 & \leq (1+\d) \f(A).
 \end{align*}
 Choose now $\d<\eps/2$, and $U\in\cu_\d$, so that 
 $U^{-1}\in\cu_{2\d}$, and $\f(UAU^*)\leq (1+\d) \f(A)< (1+\eps) 
 \f(A)$.  Replacing $A$ with $UAU^*$, and $U$ with $U^{-1}$, we obtain
 $$
 \f(A)\leq \|U^{-1}\|^2 \f(UAU^*)\leq (1+2\d)\f(UAU^*) < 
 (1+\eps)\f(UAU^*)
 $$
 and the thesis easily follows. \\
 Assume now $A\in\ca_{+}\setminus\cj_{0+}$, then 
 $\f(A)=+\infty=\f(UAU^{*})$, because otherwise $UAU^{*}\in\cj_{0+}$, 
 so that $A=U^{-1}(UAU^{*})(U^{-1})^{*}\in\cj_{0+}$, which is absurd.  
 \qed
 
 \medskip
 
 Finally we observe that, from the proof of Lemma \ref{3.2.7} the 
 following is immediately obtained
 
 \begin{Prop}\label{3.2.8} 
	 If $A\in\ca_0$ and $\|A\|\leq 1$, then $\f(ATA^*)\leq \f(T)$, for 
	 any $T\in\cj_{0+}$.
 \end{Prop}

 \section{A construction of semicontinuous traces on 
 C$^*$-algebras} \label{sec:semicontinuous}

 The purpose of this section is to show that the 
 lower-semicontinuous semifinite regularisation of the functional 
 $\f|_\ca$ of the previous section gives a trace, namely a 
 unitarily invariant weight on $\ca$.  It turns out that this 
 procedure can be applied to any weight $\t_0$, on a unital 
 C$^*$-algebra $\ca$, which is $\eps$-invariant for $\d$-unitaries of 
 a dense $^*$-subalgebra $\ca_0$.  The particular case of the 
 functional $\f|_\ca$ is treated in the next section.  First we 
 observe that, with each weight on $\ca$, namely a functional 
 $\t_0:\ca_+\to[0,\infty]$, satisfying the property $\t_0(\l 
 A+B)=\l\t_0(A)+\t_0(B)$, $\l>0$, $A,\ B\in\ca_+$, we may associate a 
 (lower-)semicontinuous weight $\t$ with the following procedure
 \begin{equation}\label{e:weight}
	 \t(A):= \sup \{ \psi(A): \psi\in\ca^*_+,\ \psi\leq\t_0 \}
 \end{equation}
 Indeed, it is known that \cite{Combes,Stratila}
 $$
 \t(A) \equiv \sup_{\psi\in\cf(\t_0)} \psi(A) 
 $$
 where $\cf(\t_0):= \{ \psi\in\ca^*_+ : \exists\ \eps>0,\ 
 (1+\eps)\psi<\t_0 \}$.  Moreover the following holds

 \begin{Thm}\label{3.3.1} {\rm \cite{QV}} 
	 The set $\cf(\t_0)$ is directed, namely, for any $\psi_1,\ 
	 \psi_2\in\cf(\t_0)$, there is $\psi\in\cf(\t_0)$, s.t. $\psi_1,\ 
	 \psi_2\leq \psi$.
 \end{Thm}

 From this theorem easily follows

 \begin{Cor}\label{3.3.2} 
	 Let $\t_0$ be a weight on the C$^*$-algebra $\ca$, and $\t$ be 
	 defined as in $(\ref{e:weight})$.  Then 
	 \itm{i} $\t$ is a semicontinuous weight on $\ca$ 
	 \itm{ii} $\t=\t_0$ iff $\t_0$ is semicontinuous.  
	 \itm{iii} The domain of $\t$ contains the domain of $\t_0$.
	 \\
	 The weight $\t$ will be called the semicontinuous regularization 
	 of $\t_0$.
 \end{Cor}
 \begin{proof}
	 $(i)$ From Theorem \ref{3.3.1}, $\t(A) = \sup_{\psi\in\cf(\t_0)} 
	 \psi(A) = \lim_{\psi\in\cf(\t_0)} \psi(A)$, whence linearity and 
	 semicontinuity of $\t$ easily follow.  \\
	 $(ii)$ is a well known result by Combes \cite{Combes}.\\
	 $(iii)$ Immediately follows from the definition of $\t$.
 \end{proof}

 \begin{Prop}\label{3.3.3} 
	 Let $\t_0$ be a weight on $\ca$ which is $\eps$-invariant by 
	 $\d$-unitaries in $\ca_0$ (as in Proposition \ref{3.2.6}).  Then 
	 the associated semicontinuous weight $\t$ satisfies the same 
	 property.
 \end{Prop}
 \begin{proof}
	 Fix $\eps<1$ and choose $\d\in(0,1/2)$, s.t. $U\in\cu_\d$ implies 
	 $|\t_0(UAU^*)-\t_0(A)|<\eps \t_0(A)$, $A\in\ca_+$.  Then, for any 
	 $U\in\cu_{\d/2}$ and any 
	 $\psi\in\ca^*_+$, $\psi\leq\t_0$, we get
	 $$
	 \psi\circ adU(A)\leq
	 \t_0(UAU^*)\leq (1+\eps)\t_0(A),
	 $$
	 for $A\in\ca_+$, $i.e.$ $(1+\eps)^{-1}\psi\circ adU\leq \t_0$.  
	 Then
	 \begin{align*}
		 \t(UAU^*) 
		 & = (1+\eps) \sup_{\psi\leq\t_0} (1+\eps)^{-1}\psi\circ adU(A) \\
		 & \leq (1+\eps) \sup_{\psi\leq\t_0} \psi(A) \\
		 & = (1+\eps)\t(A).
	 \end{align*}
	 Since $U^{-1}\in\cu_\d$, replacing $U$ with $U^{-1}$ and $A$ with 
	 $UAU^*$, we get $\t(A)\leq (1+\eps)\t(UAU^*)$.  Combining the 
	 last two inequalities, we get the result.  
 \end{proof}

 \begin{Thm}\label{3.3.4} 
	 The semicontinuous weight $\t$ of Proposition~\ref{3.3.3} is a trace on 
	 $\ca$, namely, setting $\cj_+:= \{A\in\ca_+: \t(A)<\infty\}$, and 
	 extending $\t$ to the linear span $\cj$ of $\cj_+$, we get 
	 \itm{i} $\cj$ is an ideal in $\ca$ \itm{ii} $\t(AB)=\t(BA)$, for 
	 all $A\in\cj$, $B\in\ca$.
 \end{Thm}
 \begin{proof}
	 $(i)$ Let us prove that $\cj_+$ is a unitary invariant face in 
	 $\ca_+$, and it suffices to prove that $A\in\cj_+$ implies 
	 $UAU^*\in\cj_+$, for all $U\in\cu(\ca)$, the set of unitaries in 
	 $\ca$.  Suppose on the contrary that there is $U\in\cu(\ca)$ s.t. 
	 $\t(UAU^*)=\infty$.  Then there is $\psi\in\ca^*_+$, 
	 $\psi\leq\t_0$, s.t. $\psi(UAU^*) > 2\t(A)+2$.  Then we choose 
	 $\d<3$ s.t. $V\in\cu_\d$ implies $\t(VAV^*)\leq 2\t(A)$, and an 
	 operator $U_0\in\ca_0$ s.t. $\|U-U_0\|< \min\{ \frac{\d}{3}, 
	 \frac1{3\|A\|\|\psi\|} \}$.  The inequalities
	 $$
	 \|U_0U_0^*-1\| = \|U^*U_0U_0^*-U^*\| \leq 
	 \|U^*U_0-1\|\|U_0^*\|+\|U_0^*-U^*\| < \d
	 $$
	 and analogously for $\|U_0^*U_0-1\|<\d$, show that 
	 $U_0\in\cu_\d$.  Then, since $|\psi(U_0AU_0^*)-\psi(UAU^*)|\leq 
	 3\|\psi\|\|A\|\|U-U_0\|<1$, we get
	 $$
	 2\t(A)\geq \t(U_0AU_0^*) \geq \psi(U_0AU_0^*) \geq
	 \psi(UAU^*) - 1 \geq 2\t(A)+1
	 $$
	 which is absurd. \\
	 $(ii)$ We only have to show that $\t$ is unitary invariant.  Take 
	 $A\in\cj_+$, $U\in\cu(\ca)$.  For any $\eps>0$ we may find a 
	 $\psi\in\ca^*_+$, $\psi\leq\t_0$, s.t. 
	 $\psi(UAU^*)>\t(UAU^*)-\eps$, as, by $(i)$, $\t(UAU^*)$ is 
	 finite.  Then, arguing as in the proof of $(i)$, we may find 
	 $U_0\in\ca_0$, so close to $U$ that
	 \begin{align*}
		 & |\psi(U_0AU_0^*)-\psi(UAU^*)|<\eps \\
		 & (1-\eps)\t(A)\leq \t(U_0AU_0^*) \leq (1+\eps)\t(A). 
	 \end{align*}
	 Then
	 \begin{align*}
		 \t(A) & \geq \frac1{1+\eps}\ \t(U_0AU_0^*) \geq 
		 \frac1{1+\eps}\ \psi(U_0AU_0^*) \\
		 & \geq \frac1{1+\eps}\ (\psi(UAU^*) -\eps) \geq 
		 \frac1{1+\eps}\ (\t(UAU^*) -2\eps).
	 \end{align*}
	 By the arbitrariness of $\eps$ we get $\t(A)\geq \t(UAU^*)$.  
	 Replacing $A$ with $UAU^*$, we get the thesis.  
 \end{proof}
  
 The second regularization we need turns $\t$ into a (lower 
 semicontinuous) semifinite trace, namely guarantees that
 $$
 \t(A) = \sup\{ \t(B) : 0\leq B\leq A,\ B\in\cj_+ \}
 $$
 for all $A\in\ca_+$.  In particular the semifinite regularization 
 coincides with the original trace on the domain of the latter.  This 
 regularization is well known (see $e.g.$ \cite{DixmierC}, Section 6), 
 and amounts to represent $\ca$ via the GNS representation $\pi$ 
 induced by $\t$, define a normal semifinite faithful trace $tr$ on 
 $\pi(\ca)''$, and finally pull it back on $\ca$, that is 
 $tr\circ\pi$.  It turns out that $tr\circ\pi$ is (lower 
 semicontinuous and) semifinite on $\ca$, $tr\circ\pi\leq\t$, and 
 $tr\circ\pi(A)=\t(A)$ for all $A\in\cj_+$, that is $tr\circ\pi$ is a 
 semifinite extension of $\t$, and $tr\circ\pi=\t$ iff $\t$ is 
 semifinite.

 We still denote by $\t$ its semifinite extension.  As follows from 
 the construction, semicontinuous semifinite traces are exactly those 
 of the form $tr\circ\pi$, where $\pi$ is a tracial representation, 
 and $tr$ is a n.s.f. trace on $\pi(\ca)''$.

 \section{The semicontinuous trace on  almost 
 local operators and the Novikov-Shubin numbers} \label{sec:regularized}

 Now we apply the regularization procedure described in the previous 
 section to Roe's functional.  Let us remark that the 
 semicontinuous regularization of the weight $\f|_\ca$ is a trace in 
 the sense of property $(ii)$ of Theorem~\ref{3.3.4}, which is 
 stronger then the trace property in \cite{Roe}.  First we observe 
 that $\f|_\ca$ is not semicontinuous.

 \begin{Prop}\label{3.4.1} 
	 The set $\cn_0:=\{T\in\ca_+: \f(T)=0\}$ is not closed.  In 
	 particular, there are operators $T\in\ca_+$ s.t. $\f(T)=1$ but 
	 $\t(T)=0$ for any (lower-)semicontinuous trace $\t$ dominated by 
	 $\f|_\ca$.
 \end{Prop}
 \begin{proof}
	 Recall from Lemma \ref{2.1.2}$(i)$ that there are positive real 
	 functions $\b_1,\ \b_2$ s.t. $0<\b_1(r)\leq V(x,r) \leq \b_2(r)$, 
	 for all $x\in M$, $r>0$, and $\lim_{r\to0}\b_2(r)=0$.  Therefore 
	 we can find a sequence $r_n\dec0$ s.t. $\sum_{n=1}^\infty 
	 \b_2(r_n)<\infty$.  Fix $o\in M$, and set $X_n:= \{ (x_1,x_2)\in 
	 M\times M : n\leq \d(x_i,o)\leq n+1,\ \d(x_1,x_2)\leq r_n \}$, 
	 $Y_n:= \cup_{k=1}^n X_k$, $n\leq\infty$, and finally let $T_n$ be 
	 the integral operator whose kernel, a section of $End(F)$ denoted 
	 $k_n$, is the characteristic function of $Y_n$.  Since $k_n$ has 
	 compact support, if $n<\infty$, $\f(T_n)=0$.  On the contrary, 
	 since $Y_\infty$ contains the diagonal of $M\times M$, clearly 
	 $\f(T_\infty)=1$.  Finally, by Lemma \ref{3.1.4},
	 \begin{align*}
		 \|T_\infty-T_n\| 
		 & \leq \sup_{x\in M} \sum_{k=n+1}^\infty \int_M 
		 \chi_{X_k}(x,y)dy \\
		 & \leq \sup_{x\in M} \sum_{k=n+1}^\infty V(x,r_k) 
  \leq \sum_{k=n+1}^\infty \b_2(r_k) \to0. 
	 \end{align*}
	 This proves both the assertions.
 \end{proof}
 
 \begin{Thm}\label{def:trace}
	 Let $M$ be an open manifold of C$^{\infty}$-bounded geometry, 
	 endowed with a regular exhaustion $\ck$, and $F$ a finite 
	 dimensional hermitian vector bundle on $M$.  Then the regularisation 
     $Tr_\ck$ of $\f$ is a semicontinuous semifinite trace on 
	 the C$^{*}$-algebra $\ca(F)$ of almost local operators on 
	 $L^{2}(F)$, which vanishes on compact operators.
 \end{Thm}
 \begin{proof}
	 The first statement follows by the results of the previous 
	 section.  Clearly $\f$, hence $Tr_\ck$, vanishes of finite rank 
	 operators.  By semicontinuity, the kernel of $Tr_\ck$ is closed, 
	 which implies the thesis.
 \end{proof}	 

 Finally we give a sufficient criterion for a positive operator $A\in 
 \ca$ to satisfy $Tr_\ck(A)=\f(A)$.

 \begin{Prop}\label{3.4.2} 
	 Let $A\in\cj_{0+}$ be an integral operator, whose kernel $a(x,y)$ 
	 is a section of $End(F)$ which is uniformly continuous in a 
	 neighborhood of the diagonal in $M\times M$, namely
	 \begin{equation}\label{e:uniform}
		 \forall \eps>0,\ \exists \d_\eps>0\ : \d(x,y)<\d \imply 
		 |a(x,y)-a(x,x)|<\eps.
	 \end{equation}
	 Then $Tr_\ck(A)=\f(A)$.
 \end{Prop}
 \begin{proof}
	 Consider first a family of integral operators $B_\d$, with 
	 kernels, which are sections of $End(F)$, given by
	 \begin{equation*}
		 b_\d(x,y):= \frac{\b_1(\d)}{\b_2(\d)}\ 
		 \frac{\chi_{\D_\d}(x,y)}{ V(x,\d)},
	 \end{equation*}
	 where $\D_\d := \{ (x,y)\in M\times M : \d(x,y)<\d \}$.  
	 Set $E_n$ for the multiplication operator by the characteristic 
	 function of $K_{n}$, and observe that
	 \begin{equation*}
		 Tr(E_nB_\d B_\d^*E_n)  = 
   \frac{\b_{1}(\d)^{2}}{\b_{2}(\d)^{2}}\int_{K_{n}}\frac{dx}{V(x,\d)} 
   \leq\frac{vol(K_{n})}{\b_2(\d)}
	 \end{equation*}
	 Therefore $\f(B_\d B_\d^*)\leq \b_2(\d)^{-1}$, hence 
	 $\psi_\d := \f(B_\d\cdot B_\d^*)$ is a positive functional on $\ca$. 
	 Since $\sup_{x\in M} \int_M b_\d(x,y)dy = \frac{\b_1(\d)}{ \b_2(\d)} 
	 \leq 1$, and $\sup_{y\in M} \int_M b_\d(x,y)dx \leq 
	 \frac{\sup_{y\in M} V(y,\d)}{\b_2(\d)} \leq 1$, Riesz-Thorin theorem implies 
	 $\|B_\d\|\leq1$, hence $\psi_\d\leq \f|_\ca$ by Proposition~\ref{3.2.8}.  
  By definition of $Tr_\ck$, we have $\psi_\d \leq Tr_\ck$.  \\
	 Take now $A\in\ca_+$ satisfying (\ref{e:uniform}),  
	 $\eps>0$, $\d=\d_\eps<\eps$, and set $\b(\d):= (\frac{\b_1(\d) }{ 
	 \b_2(\d)})^2$. Then
	 \begin{align*}
		 &|Tr(E_nB_\d AB_\d^*E_n) - \b(\d) Tr(E_nAE_n)|\\
		 &\leq \int_{K_{n}} dx \int_{B(x,\d)\times B(x,\d)} 
		 b_\d(x,y)|a(y,z)-a(x,x)|b_\d(x,z) dydz \\
		 &\leq 3\eps \int_{K_{n}} dx \int_{B(x,\d)\times B(x,\d)} 
		 b_\d(x,y)b_\d(x,z) dydz \leq 3\eps \b(\d) vol(K_{n}),
	 \end{align*}
	 hence $|\f(A)-Tr_\ck(A)|\leq|\f(A)-\psi_\d(A)| 
 \leq 3\eps \b(\d) + (1-\b(\d))\f(A)$.  The thesis follows by 
	 Lemma \ref{2.1.2}$(ii)$.
 \end{proof}

 \begin{Thm}\label{3.4.3} 
	 Let $M$ be an open manifold of C$^{\infty}$-bounded geometry, 
	 endowed with a regular exhaustion $\ck$, $p\in\{0,1,\ldots,dim 
	 M\}$ and $Tr_\ck$ the trace on $\ca(\La^{p}T^{*}M)$ given by 
	 Theorem \ref{def:trace}.  Then $\e{-t\D_p}$ belongs to the 
  domain of $Tr_\ck$, for any $t>0$, and
	 $$
	 Tr_\ck(\e{-t\D_{p}}) = \lo \frac{\int_{K_{n}} tr(H_{p}(t,x,x)) 
	 dx}{vol(K_{n})},
	 $$
	 where $H_{p}$ is the heat kernel on $p$-forms.
 \end{Thm}
 \begin{proof}
	 By Remark~\ref{rem:laplacian} and Theorem~\ref{3.1.6} we have 
	 that $\e{-t\D_p}$ belongs to $\cj_{0+}\cap\ca$, hence, by 
	 Corollary~\ref{3.3.2}, $(iii)$, it belongs to the domain of 
	 $Tr_\ck$.  The equality then follows by Propositions \ref{2.1.5}, 
	 \ref{3.4.2}, and \ref{p:Roeformula}. \\
	  \end{proof}

 We note that the previous result can be easily extended to functional calculi
 $f(\Delta)$ where $f$ has exponential decay.
 
 \medskip
 
 We conclude this section showing that the above defined trace can be 
 used to define the Novikov-Shubin numbers for open manifolds.  Let 
 $M$ be an open $n$-manifold with C$^{\infty}$-bounded geometry 
 possessing a regular exhaustion, and denote by 
 $\ca_p\equiv\ca(\La^{p}T^{*}M)$ the C$^{*}$-algebra of almost local 
 operators acting on $L^{2}(\La^{p}T^{*}M)$.  Let us define a positive 
 measure $\m_{p}$ on $[0,\infty)$ as
 $$
 Tr_\ck(f(\Delta_{p}))=\int f(t) d\m_{p}(t),\quad f\in C_{c}[0,\infty),
 $$
 and denote by $N_{p}(t) := \m_{p}([0,t))$ the spectral density function 
 of $\Delta$, and by $\th_{p}$ the function 
 $\th_p(t):=Tr_\ck(e^{-t\Delta_{p}})=\int_{0}^{\infty} \e{-t\l} d\m_{p}(\l)$.  
 Then $\th_{p}$ can be written as $\th_p(t) = \int_{0}^{\infty}\e{-t\l} 
 dN_p(\l)$ so that, by a Tauberian theorem (\cite{GS}, Appendix), 
 $\lim_{t\to0} N_p(t)= \lim_{t\to\infty}\th_p(t)$.
  
 \begin{Dfn} \label{def:NS-inv} 
	 We define $b_{p}\equiv b_p(M,\ck) :=\lim_{t\to0} N_{p}(t)= 
	 \lim_{t\to\infty}\th_{p}(t)$ to be the $p$-th L$^{2}$-Betti 
	 number of the open manifold $M$ endowed with the exhaustion 
	 $\ck$.  Let us now set $N^{0}_{p}(t) := N_{p}(t) - b_{p} \equiv 
	 \m_{p}((0,t))$, and 
	 $\th^{0}_{p}(t) := \th_{p}(t) - b_{p} = \int_{0}^{\infty}\e{-t\l} 
	 dN^{0}_{p}(\l)$.  The Novikov-Shubin numbers of $(M,\ck)$ are 
	 then defined as
	 \begin{align*}
		 \a_{p}\equiv \a_{p}(M,\ck) & := 2\limsup_{t\to0} \frac{\log 
		 N^{0}_{p}(t)}{\log t}, \\
		 \underline{\a}_{p}\equiv \underline{\a}_{p}(M,\ck) & := 
		 2\liminf_{t\to0} \frac{\log N^{0}_{p}(t)}{\log t}, \\
		 \a'_{p}\equiv \a'_{p}(M,\ck) & := 2\limsup_{t\to\infty} 
		 \frac{\log \th^{0}_{p}(t)}{\log 1/t}, \\
		 \underline{\a}'_{p}\equiv \underline{\a}'_{p}(M,\ck) & := 
		 2\liminf_{t\to\infty} \frac{\log \th^{0}_{p}(t)}{\log 1/t}.
	 \end{align*}
 \end{Dfn}
 
 It follows from (\cite{GS}, Appendix) that $\underline{\a}_{p} = 
 \underline{\a}'_{p}\leq \a'_{p}\leq \a_{p}$, and 
 $\a'_{p} = \a_{p}$ if $\th_{p}^{0}(t)= O(t^{-\d})$, for 
 $t\to\infty$, or equivalently $N^{0}_{p}(t)=O(t^{\d})$, for $t\to0$.
 Observe that L$^{2}$-Betti numbers and Novikov-Shubin numbers depend 
 on the limit procedure $\om$ involved in the trace, and on the exhaustion
$\ck$.

 \begin{rem} $(i)$ $L^2$-Betti numbers for open manifolds have been 
     defined by J.~Roe (\cite{Roe}) as $\inf \{ 
     \phi(f(\D_{p})) : f\in C_{c}^{\infty}[0,\infty)_{+}, f(0)=1 \}$, who
     showed in \cite{RoeBetti} their invariance under quasi-isometries.
     One can show easily that the two definitions coincide. 
     \\
     $(ii)$  In \cite{GI4} we used the semicontinuous semifinite 
     trace constructed here to show that Novikov-Shubin numbers are 
     asymptotic dimensions, and showed that they are invariant under
     quasi-isometries.  Observe that the definition of the 
     Novikov-Shubin numbers used there makes use of the noncommutative 
     Riemann integration, which allows us to take the trace of many 
     projections.  Since the spectral projection 
     $\chi_{[0,t)}(\Delta_{p})$ is Riemann integrable for almost all 
     $t$, the two definitions coincide.
     \\
     $(iii)$ Corollary 3.16 in \cite{GI5} shows that the above definitions 
     for $L^2$-Betti numbers and Novikov-Shubin numbers coincide with the
     classical ones in the case of amenable coverings, if one chooses the 
	    exhaustion given by the F\o lner condition.
     \end{rem}
     
     \medskip
  
  In the case of coverings there is a well-known conjecture 
	 on the positivity of the $\a_{p}$'s.  A result by Varopoulos 
	 \cite{Varopoulos1} shows that $\a_{0}$ is a positive integer, 
	 hence $\a_{0}\geq1$.  The following proposition extends this 
	 inequality to the case of open manifolds.  Moreover $\a_{0}$
  coincide with asymptotic dimension of $M$, cf. \cite{GI6}, for a suitable
  class of open manifolds, hence it can assume any value in $[1,\infty)$.
  
 \begin{Prop}\label{p:a0=a'0}
	Let $M$ be a complete non-compact Riemannian manifold of positive 
	injectivity radius and Ricci curvature bounded from below.  Then 
	$\a_{0}(M,\ck)=\a'_{0}(M,\ck)\geq 1$ for any regular exhaustion 
	$\ck$.
 \end{Prop}
 \begin{proof}
	 Recall that, under the previous assumptions, Varopoulos 
	 \cite{Varopoulos2} proved that the heat kernel on the diagonal 
	 has a uniform inverse-polynomial bound, more precisely, in the 
	 strongest form due to \cite{ChavelFeldman}, we have
	 $$
	 \sup_{x,y\in M} H_{0}(t,x,y)\leq Ct^{-1/2}
	 $$
	 for a suitable constant $C$.  Then, as 
	 $$
	 \th(t)=\t(\e{-t\D})= \lo 
	 \frac{\int_{B(o,n_{k})}H_{0}(t,x,x)dvol(x)}{V(o,n_{k})} \leq 
	 Ct^{-1/2},
	 $$
	 it follows from (\cite{GS}, Appendix) that 
	 $\a_{0}=\a_{0}'$, which concludes the proof.
 \end{proof}


\end{document}